\documentstyle[12pt]{article}
   \font\tenmsb=msbm10 scaled\magstep 1
   \font\sevenmsb=msbm7 scaled \magstep 1
   \font\fivemsb=msbm5 scaled \magstep 1
\newfam\msbfam
      \textfont\msbfam=\tenmsb
      \scriptfont\msbfam=\sevenmsb
      \scriptscriptfont\msbfam=\fivemsb
\def\Bbb#1{{\fam\msbfam #1}}

\font\tengothic=eufm10 scaled\magstep 1
\font\sevengothic=eufm7 scaled\magstep 1
\newfam\gothicfam
      \textfont\gothicfam=\tengothic
      \scriptfont\gothicfam=\sevengothic

\newcommand{\verylong}{\hbox{$\hbox to .45in{\rightarrowfill}$}  }
\newcommand{\kindalong}{\hbox{$\hbox to .35in{\rightarrowfill}$}  }

\newtheorem{thm}{Theorem}[section]
\newtheorem{prop}[thm]{Proposition}
\newtheorem{lemma}[thm]{Lemma}
\newtheorem{cor}[thm]{Corollary}

\newtheorem{defn}[thm]{Definition}
\newtheorem{rmk}[thm]{Remark}
\newtheorem{exemp}[thm]{Example}

\newcommand{\qed}{\mbox{\hskip 1cm $\rlap{$\sqcap$}\sqcup$}}

\newenvironment{proof}{{\em Proof: }}{\qed\bigskip}
\newenvironment{example}{\begin{exemp}\em}{\end{exemp}}
\newenvironment{remark}{\begin{rmk}\em}{\end{rmk}}
\newenvironment{definition}{\begin{defn}\rm}{\end{defn}}

\newcommand{\diam}{\mathop{\rm diam\,}}
\newcommand{\reg}{\mathop{\rm reg\,}}
\newcommand{\pthree}{{\Bbb P}^3}
\newcommand{\pn}{{\Bbb P}^n}
\newcommand{\IC}{{\cal I}_C}
\newcommand{\hover}{\overline{h}}

\title{Degrees of Generators of Ideals Defining \\
	Subschemes of Projective Space}
\author{  Heath M. Martin \\
	  Florida State University \\
	  Tallahassee, FL  32306 \\
	  {\normalsize{\em e-mail}: martin@math.fsu.edu}
	\and
	  Juan C. Migliore \\
	  University of  Notre Dame \\
	  Notre Dame, IN  46556 \\
	  {\normalsize{\em e-mail}: Juan.C.Migliore.1@nd.edu}
 	}
\date{}

\begin{document}
\setlength{\baselineskip}{18pt}
\maketitle

In this paper, we are interested in giving an upper bound for the
degrees of the generators for an ideal defining a curve in projective
space, and in investigating the properties of curves for which
this bound is realized.  While our bound works for subschemes
of any dimension in any projective space,
our characterizations for when the bound is realized
only work for curves in $\pthree$.  Using the known structure
theory for liaison classes of curves in $\pthree$ allows us to
give a reasonably complete picture of the liaison classes containing
the extremal curves, mainly in terms of the dimensions of the components
of the deficiency module for the liaison class.  We are also able
to give some similar conditions for when a liaison class contains
an integral curve satisfying our upper bound for degrees of generators.

In contrast to the minimal degree of a generator, in most cases,
the maximal degree is not readily computed in terms of other invariants,
for example, the Hilbert function.  There is, of course,
Castelnuovo--Mumford regularity, \cite{mumford}, bounding
the degrees of generators in terms of the cohomology modules
(and more generally, bounding the degrees of generators of
the syzygies). Another result which gives
information about a certain class of curves is in \cite{DGM}, where they
show that the defining ideals of arithmetically Cohen--Macaulay curves
have generators whose degree is bounded above by one plus the last
integer for which the second difference of the Hilbert function
is non-zero.  Furthermore, in \cite{CGO}, most arithmetically Cohen--Macaulay
curves achieving this bound are classified;
generally speaking, with some exceptions,  the arithmetically
Cohen--Macaulay curves whose ideals have a generator of maximal degree
are linked in the minimal number of steps to  plane curves.  In that
paper, however, they note that little is known about the non-arithmetically
Cohen--Macaulay case, or the non-codimension two case.  This is what
sparked our interest in the problem, and we offer a relatively complete
solution in the present paper, at least for the case of curves.

Now, we describe the contents of this paper more precisely.  In the
first section, we set up our notation and state and prove the
upper bound for curves in $\pn$ (Proposition~\ref{maxdegree}),
and more generally for  subschemes of arbitrary codimension in $\pn$
(Proposition~\ref{general_maxdegree}).  We also interpret
the bound as a condition on the cohomology of the curve,
which we designate {\it equal cohomology}.  In the
second section, we specialize to $\pthree$, and study the
even liaison classes of curves in an effort to determine
those classes having a curve with the equal cohomology property
necessary for having a high degree generator.  This classification
depends heavily on the structure theory for even liaison classes.
Our main result of this section is Theorem~\ref{coh-prop},
which we use to characterize the liaison classes having
curves with equal cohomology, in terms of the Hilbert
function of the minimal curve in the class.
In the third section, we show that the curves with equal
cohomology enjoy a strong Lazarsfeld--Rao structure, in the
sense that liaison classes of such curves have minimal curves,
and the other curves are obtained by basic double links of
low degree and deformations.  The fourth section is devoted to
describing which of the curves with equal cohomology can possibly
be integral, by using the machinery developed by Nollet of
examining the postulation characters of curves.  Finally, in
the last section, we reconsider the problem of generators
of high degree, and interpret our previous results in
this context.  As a result, we are able to give a
mostly complete picture in terms of cohomology of which
curves have high degree generators.  This picture is
particularly clear for the Buchsbaum curves, since their
cohomology is quite well-known.  We are also able to
give some necessary conditions for integral curves to
have generators of high degree, and again, the Buchsbaum
case provides the best statement.

\section{Degrees of Generators}
Throughout this paper, we will work over an algebraically closed field $k$,
of arbitrary characteristic.  Let $S = k[x_0, \dots, x_n]$
be the polynomial ring over $k$, and $\pn$ the $n$-dimensional
projective space over $k$.  We will furthermore consider only
subschemes of $\pn$ which are locally Cohen--Macaulay and equidimensional.
It is well-known that this is equivalent to requiring that all
the intermediate cohomology modules $H^i_*(\pn, V)$, $1 \leq i \leq \dim V$,
have finite length.

Given a subscheme $V$ of $\pn$, with $\dim V = d$, let $I_V$ denote
its homogeneous, saturated defining ideal in $S$.  Thus $S/I_V$ is
a standard graded $k$-algebra, and so we can define the Hilbert
function of $S/I_V$ by
$$
H(S/I_V, t) = \dim_k [S/I_V]_t.
$$
Alternatively, we sometimes write $H(V, t)$ for $H(S/I_V, t)$.  It
is a standard fact that there is a polynomial $P(S/I_V, t)$, having
degree $d$, such that $H(S/I_V, t) = P(S/I_V, t)$ for
all $t \gg 0$.  We furthermore define the $n$th difference
of $H(S/I_V, t)$ inductively as follows:
\begin{eqnarray*}
\Delta^1 H(S/I_V,t) &=& H(S/I_V, t) - H(S/I_V, t-1) \\
\Delta^n H(S/I_V,t) &=& \Delta^1 (\Delta^{n-1} H(S/I_V, t)).
\end{eqnarray*}
Now, since $H(S/I_V, t)$ is eventually a polynomial of degree $d$,
the function $\Delta^{d+1} H(S/I_V, t)$ is eventually zero, and we
define
$$
\sigma(S/I_V) = \min \{\, k : \Delta^{d+1} H(S/I_V, t) = 0 \mbox{ for all
				$t \geq k$} \, \}.
$$
Again, we will sometimes write $\sigma(V)$ for $\sigma(S/I_V)$.
It is worth noting that if $\dim V = d$, then
the Hilbert function of $V$ and the Hilbert polynomial
of $V$ are equal in all degrees $\geq t + d + 1$ if
and only if $\sigma(V) = t$.
The degree at which the Hilbert function and
the Hilbert polynomial agree from then on is sometimes called,
at least in the context of local algebra, the postulation number,
and has played an important role in
questions about Cohen--Macaulayness and related invariants in local rings.

Given an ideal $I$, we will write $\alpha(I)$ for the minimal degree
of a minimal generator, and $\omega(I)$ for the maximal
degree of a minimal generator of $I$.

Next, let $H^i_*(\pn, {\cal I}_V) = \oplus_{s\in {\Bbb Z}} H^i(\pn, {\cal
I}_V(s))$
be the cohomology modules of $V$.  We will put
$h^i(\pn, {\cal I}_V(t)) = \dim_k H^i(\pn, {\cal I}_V(t))$.
We let $e(V) = \max\{\, t : h^{d+1}(\pn, {\cal I}_V(t)) \not=0\,\}$
denote the index of speciality of $V$.  Also, note
that by our assumption that subschemes be locally Cohen--Macaulay
and equidimensional, $h^i(\pn, {\cal I}_V(t))$ is non-zero
for only finitely many $t$, when $1 \leq i \leq d$.
Hence for a non-arithmetically Cohen--Macaulay curve $C$ in $\pn$,
with notation following Martin-Deschamps and Perrin, we can define
$$
r_a(C) = \min \{\, n \in {\Bbb Z} : h^1(\pn, \IC(n)) \not = 0 \,\} \quad\quad
r_o(C) = \max \{\, n \in {\Bbb Z} : h^1(\pn, \IC(n)) \not = 0 \,\},
$$
and $\diam H^1_*(\pn, C) = r_o(C) - r_a(C) + 1$, the number of
components between the first and the last non-zero components, inclusive.
Note that some of the intermediate components may have dimension zero,
but we also count them.  If $C$ is arithmetically Cohen--Macaulay,
then $H^1_*(\pthree, C) = 0$, and we will put $\diam H^1_*(\pthree, C) = 0$.

In this section, we prove a statement about the maximal degree
of a generator for the defining ideal of a curve,
in terms of the Hilbert function and the cohomology
of the curve.  The relationship between these two objects is well-known,
and we spell it out explicitly in the first lemma.

\begin{lemma} \label{hilbfcn_eqn}
Let $C$ be a curve in $\pn$.  Then
\begin{eqnarray*}
\Delta^2 H(C,t) &=& h^2(\pn, \IC(t)) - 2h^2(\pn, \IC(t-1)) + h^2(\pn, \IC(t-2))
\\
&& \quad\quad \mbox{} - h^1(\pn, \IC(t)) + 2h^1(\pn, \IC(t-1)) - h^1(\pn,
\IC(t-2))
\end{eqnarray*}
\end{lemma}

\begin{proof}  First, recall that the Hilbert polynomial is given
by $P(C,t) = h^0(\pn, {\cal O}_C(t)) - h^1(\pn, {\cal O}_C(t))$
(see \cite[Exercise III.5.2]{hart}).  Thus,
from  the short exact sequence
$$
0 \rightarrow I_t \rightarrow S_t \rightarrow H^0(\pn, {\cal O}_C(t))
	\rightarrow H^1(\pn, \IC(t)) \rightarrow 0,
$$
we obtain
\begin{eqnarray*}
H(C,t) &=& h^0(\pn, {\cal O}_C(t)) - h^1(\pn, \IC(t)) \\
	&=& P(C,t) + h^2(\pn, \IC(t)) - h^1(\pn, \IC(t)).
\end{eqnarray*}
Now, since $P(C, t)$ is a polynomial of degree $1$, on taking second
differences we obtain
\begin{eqnarray*}
\Delta^2 H(C,t) &=& h^2(\pn, \IC(t)) - 2h^2(\pn, \IC(t-1)) + h^2(\pn, \IC(t-2))
\\
&& \quad\quad \mbox{} - h^1(\pn, \IC(t)) + 2h^1(\pn, \IC(t-1)) - h^1(\pn,
\IC(t-2))
\end{eqnarray*}
which is what we wanted to show.
\end{proof}

The following corollary is an immediate consequence.

\begin{cor}  \label{Delta-coh}
$\Delta^2 H(C, t) = 0$ for all $t \geq k+2$ if and only
if $h^2(\pn, \IC(t)) = h^1(\pn, \IC(t))$ for all $t \geq k$. \qed
\end{cor}

\begin{prop} \label{equal-coh}
Suppose $C$ is a curve in $\pn$ defined by an ideal $I = I_C$.
If $\omega(I) = \sigma(C) + k$ for some $k \geq 1$, then
$e(C) = r_o(C)$ and $h^2(\pn, \IC(t)) = h^1(\pn, \IC(t))$
for $t \geq e(C) -k+1$.
\end{prop}

\begin{proof}  We first show that $e(C) = r_o(C)$.  If not,
let $m = \max \{\, e(C), r_o(C) \,\}$.  Then by Castelnuovo--Mumford
regularity \cite{mumford}, we have
$$
\sigma(C) < \sigma(C) + k = \omega(I) \leq \reg(I) \leq m+3.
$$
But $\Delta^2 H(C, m+2) \not=0$, because of Lemma~\ref{hilbfcn_eqn}
and so $\sigma(C) = m+3$, which is a contradiction.
Thus, we must have $e(C) = r_o(C) = m$.  In particular,
$\reg(I) = m+3$, and so again by Castelnuovo--Mumford
regularity, we have $\sigma(C) = \omega(I) - k \leq \reg(I) - k = m + 3 - k$.
Thus, we have $\Delta^2 H(C, t) = 0$ for $t \geq m+3-k$, and so by
Lemma~\ref{hilbfcn_eqn}, $h^1(\pn, \IC(t)) = h^2(\pn, \IC(t))$
for $t \geq m-k+1$.
\end{proof}

\begin{prop}  \label{maxdegree}
Suppose $I = I_C$ defines a curve $C$ in ${\Bbb P}^n$.
Then $\omega(I) \leq \sigma(S/I) + \diam H^1_*(\pn, \IC)$.
\end{prop}

\begin{proof}
This follows immediately from the previous proposition, since
$h^1(\pn, \IC(t))$ and $h^2(\pn, \IC(t))$ can only possibly be
non-zero and
equal for $k=\diam H^1_*(\pn, \IC)$ degrees.
\end{proof}

\begin{remark}  Also note that Proposition~\ref{maxdegree} includes
the case that $C$ is arithmetically Cohen--Macaulay, and says
that $\omega(I_C) \leq \sigma(S/I_C)$.  This is a result in
\cite{DGM};  see also \cite[Proposition 1.2]{CGO}.
\end{remark}

\begin{remark}  Some comments about the proof of this result are in order.
First of all, if either $e(C) \not= r_o(C)$ or if $e(C) = r_o(C)$ and
$h^1(\pn, {\cal I}_C(r_o)) \not= h^2(\pn, {\cal I}_C(r_o))$, we get
$\omega(I) \leq \sigma(S/I)$.  This is the same bound as when $C$
is assumed to be arithmetically Cohen--Macaulay.  Thus,
the cases of most interest occur when $H^1_*(\pn, {\cal I}_C)$ and
$H^2_*(\pn, {\cal I}_C)$ both become zero at the same degree, and
moreover have equal dimensions for some number of preceding degrees.
Essentially, Proposition~\ref{equal-coh} says that
having a generator of high degree forces $h^1(\pn, {\cal I}_C(t))$
and $h^2(\pn, {\cal I}_C(t))$ to be equal in a large number of degrees.
\end{remark}

Curves which are not arithmetically Cohen--Macaulay and have a generator
of maximum degree in the sense of the above proposition, must be
``almost Buchsbaum.''  This means that a general linear form $L$
induces a multiplication on $H^1_*(\pn, \IC)$ which has non-trivial
kernel in each degree.  As notation, if $L$ is a linear form,
let $K_L$ be the kernel of the multiplication on $H^1_*(\pn, \IC)$
induced by $L$.

\begin{prop}  \label{almostBuchsbaum}
Let $C$ be a non-arithmetically Cohen--Macaulay curve
in $\pthree$, and suppose that $h^1(\pn, \IC(t)) = h^2(\pn, \IC(t))$
in the last $r$ degrees.  Let $L$ be a general linear form and
let $K = K_L$.  Then $\dim_k K_{t} > \dim_k K_{t+1} > 0$ for
all $t = r_o(C) - r + 1, \dots, r_o(C) - 1$.
\end{prop}

\begin{proof}  Let $L$ be a general general linear form defining
a hyperplane $H$, and let $Z = C \cap H$ be the hyperplane section
of $C$ considered as a subscheme of ${\Bbb P}^{n-1}$.  Then it is easy to
see that
$$
\Delta^2 H(C, t) = \Delta^1 H(Z, t) + \Delta^1 \dim_k K_{t-1}.
$$
Let $p = r_o(C)$, and note that by the condition on cohomology
and by Lemma~\ref{hilbfcn_eqn}, we have
$$
0 = \Delta^2 H(C, p+2) = \Delta^1 H(Z, p+2) + \Delta^1 K_{p+1}.
$$
But $\dim_k K_{p+1} = 0$ and $\dim_k K_p > 0$, so
$\Delta^1 K_{p+1} < 0$.  This implies that $\Delta^1 H(Z, p+2) > 0$,
and since $\Delta^1 H(Z, t)$ is non-increasing in the range
in which we are interested (see \cite{DGM}),
then $\Delta^1 H(Z,t) > 0$ for all
$t = r_o(C) - r + 3, \dots, r_o(C) + 2$.  But the assumptions
on the cohomology of $C$ then imply
$$
0 = \Delta^2 H(C, t) = \Delta^1 H(Z,t) + \Delta^1 \dim _k K_{t-1}
$$
for $t = r_o(C) - r + 3, \dots, r_o(C) + 2$, and so
$\Delta^1 \dim_k K_{t-1} < 0$.  That is, $\dim_k K_{t} > \dim_k K_{t+1} > 0$
for $t = r_o(C) - r+1, \dots, r_o(C) - 1$.
\end{proof}

As an immediate corollary, we obtain the following statement:

\begin{cor}  \label{almostBuchs-deg}
Suppose $C$ is a curve in $\pn$ having a generator
of degree $\sigma(C) + \diam H^1_*(\pn, \IC)$.  Then
for each general linear form $L$, $\dim K_t > \dim K_{t+1} > 0$,
for $t = r_a(C), \dots, r_o(C) - 1$. \qed
\end{cor}

We can also use Proposition~\ref{almostBuchsbaum} to refine
the bound on the maximum degree of a generator.

\begin{cor}  If $C$ is a curve in $\pn$ defined by an ideal $I = I_C$,
then $\omega(I) \leq \sigma(C) + \diam K$.
\end{cor}

\begin{proof}  Again, this follows from Proposition~\ref{equal-coh} and
Proposition~\ref{almostBuchsbaum}.
\end{proof}

More generally, we have the following result for subschemes of dimension
$d$ in $\pn$:

\begin{prop}\label{general_maxdegree}
Let $V$ be a subscheme of $\pn$ having dimension $d$,
defined by an ideal $I = I_V$.  Then
$\omega(I) \leq \sigma(V) + \max \{\, \diam H^i_*(\pn, {\cal I}_V) :
						i = 1, \dots, d\,\}$.
\end{prop}

\begin{proof} Since the proof of this result is quite similar to
the case of curves, we will only give an outline.
The Hilbert polynomial of $V$ is given by
$$
P(V, t) = \sum_{i=0}^d (-1)^i h^i(\pn, {\cal O}_V(t)),
$$
and so from the exact sequence
$$
0 \rightarrow I_t \rightarrow S_t \rightarrow H^0(\pn, {\cal O}_V(t))
	\rightarrow H^1(\pn, {\cal I}_V(t)) \rightarrow 0,
$$
we see that the Hilbert function of $V$ is
$$
H(V, t) = P(V, t) + \sum_{i=1}^d (-1)^i h^i(\pn, {\cal I}_V(t)).
$$
Since $P(V,t)$ is a polynomial of degree $d$, when we take
$(d+1)$th differences, we get
$$
\Delta^{d+1} H(V, t) = \sum_{i=1}^d (-1)^i
		\sum_{j=0}^{d+1} (-1)^j{{d+1} \choose j} h^i(\pn, {\cal I}_V(t-j)).
$$
Now we argue by cases.  If none of the cohomology modules
end in the same place, it is easy to see by Castelnuovo--Mumford
regularity that $\omega(I) \leq \sigma(V)$.
Suppose, on the other hand, that some of the cohomology modules
end in the same degree $t$, say, and the others end
in degrees $< t$, and let $m$ be the
maximum of the diameters of the intermediate cohomologies.
Then $\reg(V) = t + d + 1$, and by the formula above,
$\sigma(V) \geq t - m +d + 1$, and since $\omega(I) \leq \reg(I)$,
the required inequality follows.
\end{proof}

We can be a bit more precise in a few cases.  For instance,
suppose $V$ is a surface in $\pn$.  Then there are only
three non-zero cohomology modules, and as in the proof above
the Hilbert function of $V$ is given by
$$
H(V, t) = P(V, t) - h^3(\pn, {\cal I}_V(t)) + h^2(\pn, {\cal I}_V(t))
				- h^1(\pn, {\cal I}_V(t)).
$$
Because the $h^3$ term and the $h^1$ term have the same sign,
the only cancellation that can occur comes from the $h^2$ term, and using
the same argument as above, we get the inequality
$\omega(I) \leq \sigma(I) + \diam H^2_*(\pn, {\cal I}_V)$.
In particular, if $V$ were a non-arithmetically Cohen--Macaulay
surface, with $H^2_*(\pn, {\cal I}_V) = 0$, then
$\omega(I) \leq \sigma(V)$, which is the same bound
as in the arithmetically Cohen--Macaulay case.

\section{Curves with equal cohomology}

The previous section showed that the property of having a generator
of high degree is very closely related to having equal (non-zero) cohomology
dimensionally in a large number of degrees.  This section is devoted
to characterizing when an even liaison class of curves in $\pthree$
has this property.  By a slight abuse of terminology, we make the
following definition:

\begin{defn}  A curve $C$ in $\pthree$ is said to have {\em equal
cohomology} if $e(C) = r_o(C)$ and $h^1(\pthree, \IC(t)) = h^2(\pthree,
\IC(t))$
for $t = r_a(C), \dots, r_o(C)$.
\end{defn}

We first recall the structure theory for curves in ${\Bbb P}^3$
(which holds more generally for codimension $2$ subschemes of $\pn$) initiated
by Lazarsfeld and Rao, and developed
in a series of papers, of which \cite{BBM} contains
the most general statement and proof.
See also the book \cite{juan-book} for a comprehensive overview
of liaison theory and the Lazarsfeld--Rao structure theory. First,
let $C$ be a curve, and choose a form $F \in I_C$ of degree $f$ and
a form $G \in S$ of degree $g$ so that $F$ and $G$ have no common
components.  Then $I_Z = G \cdot I_C + (F)$ defines a curve $Z$ in ${\Bbb
P}^3$,
called a basic double link of $C$, and denoted
$$
C:(g,f) \rightarrow Z.
$$
It is easy to see that $Z$ is evenly linked to $C$, and that there is
a short exact sequence
$$
0 \rightarrow S(-g-f) \stackrel{[F\;G]}{\longrightarrow}
	I_C(-g) \oplus S(-f) \stackrel{\phi}{\longrightarrow}
	I_Z \rightarrow 0,
$$
where $\phi(r,s) = rG+sF$.

The Lazarsfeld--Rao property says essentially that
even liaison classes of curves are built up by this process
of basic double linkage.  More precisely, let ${\cal L}$
be an even liaison class of curves in ${\Bbb P}^3$.  Then the cohomology
module $M = H^1_*(C)$ is invariant up to shifts in grading as $C$ varies
in ${\cal L}$, and there is a leftmost shift of $M$ which
actually occurs as the deficiency module of a curve in ${\cal L}$, and
every rightward shift is realized.  Thus ${\cal L}$ is parameterized
by shifts of $M$, and a curve $C_0$ which has $H^1_*(C_0)$
in the leftmost shift is called a {\em minimal curve}.  Every other
curve $C$ in ${\cal L}$ is obtained from $C_0$ as follows:  there is
a curve $C_m$ which is a deformation of $C$ through curves having
constant cohomology, and a series of basic double links
\begin{equation}\label{bdl-chain}
C_0 : (1, d_0) \rightarrow C_1 : (1, d_1) \rightarrow \cdots
	\rightarrow C_{m-1} : (1, d_{m-1}) \rightarrow C_m.
\end{equation}
We can moreover choose the degrees to satisfy
$d_0 = \cdots = d_s < d_{s+1} < \cdots < d_{m-1}$.  Note furthermore
that $C_i$ is in the $i$th shift of ${\cal L}$; that is, $H^1_*(C_i)$
is a rightward shift by $i$ degrees of $H^1_*(C_0)$.  Also,
for curves in $\pthree$, the book \cite{mdp1} gives much more
information about the behavior of invariants along liaison classes,
and moreover gives an algorithm for computing the minimal curve
in the liaison class from the deficiency module
associated to the class.

We will need to have some information about how Hilbert functions
change as we move along an even liaison class by basic double linkage.
The following result is quite elementary.

\begin{lemma} \label{hilbfcn-bdl}
Suppose $I = I_C$ defines a curve $C$ in $\pthree$.  Let $L$ be a
general hyperplane and let $F$ define a surface of degree $d$ containing
$C$.  Form the basic double link $Z$ of $C$ by $L$ and $F$.  Then
$$
\Delta^2 H(Z, t) = \left\{ \begin{array}{ll}
		\Delta^2 H(C, t-1) + 1 & \mbox{\rm if $1 \leq t \leq d-1$} \\
		\Delta^2 H(C, t-1)     & \mbox{\rm if $t \geq d$.}
		\end{array}\right.
$$
\end{lemma}

\begin{proof}  Since $Z$ is a basic double link of $C$, we have a short
exact sequence
$$
0 \rightarrow S(-d-1)  \rightarrow I_C(-1) \oplus S(-d) \rightarrow I_Z
	\rightarrow 0.
$$
Using the additive properties of Hilbert functions, it is easy to see that
$\Delta^2 H(Z, t) = \Delta^2 H(C, t-1) + \Delta^3 H(F, t)$.  Then
the statement follows, since $\Delta^3 H(F, t) = 1$ for
$0 \leq t \leq d-1$, and is zero otherwise.

Note that this is also in \cite[Corollary 2.3.5]{nollet}, in terms
of postulation characters.
\end{proof}

\begin{definition}  \label{deltadef}
Suppose there is a chain
$$
{\cal C}  : C_0 : (1, d_0) \rightarrow C_1 : (1, d_1) \rightarrow \cdots
	\rightarrow C_{m-1} : (1, d_{m-1}) \rightarrow C_m
$$
of basic double links by surfaces $F_i$ having degrees $d_i$.  Then
define $\delta({\cal C}, t, s)$ to be the number of $F_i$ such that
$d_i \geq t-s+i+1$, for $i = 1, \dots, s-1$.
\end{definition}

It is easy to see that
$\delta({\cal C}, t, s) \leq \delta({\cal C}, t-1, s)$.

The following result is a straightforward calculation using
Lemma~\ref{hilbfcn-bdl}.

\begin{cor} \label{hilbfcn-chain}
Let ${\cal C}$ as above be a chain of basic double links.  Then for
each $s = 0, \dots, m$,
$$
\Delta^2 H(C_s, t) = \Delta^2 H(C_0, t-s) + \delta({\cal C},t,s).\qed
$$
\end{cor}

The main result of this section concerns when an even liaison class has a curve
with $h^1$ and $h^2$ equal for some number of places, and
gives a characterization
in terms of the Hilbert function of the minimal curve in the liaison
class.  We first state a more general version, from which we
can trivially make a statement about minimal curves.

\begin{thm}  \label{coh-prop}  Let ${\cal L}$ be an even liaison
and let $C_0$ be a curve in ${\cal L}$.
Then a sequence
$$
C_0 \rightarrow \dots \rightarrow C_m = C
$$
of basic double links can be constructed
with $C$ having $h^1(\pthree, \IC(t)) = h^2(\pthree, \IC(t))$ in the last
$r$ places if and only if $C_0$  has
$e(C_0) \leq r_o(C_0)$ and
has Hilbert function satisfying
$$
\Delta^2 H(C_0, r_o(C_0) - r + 3) \leq \Delta^2 H(C_0, r_o(C_0) - r + 4) \leq
 	\cdots \leq \Delta^2 H(C_0, r_o(C_0) + 2) \leq 0.
$$
\end{thm}

\begin{proof}
Suppose that $C_0$ has $e(C_0) \leq r_o(C_0)$ and Hilbert function
$$
\Delta^2 H(C_0, t) = \quad \cdots \quad t_1 \quad t_2 \quad \cdots \quad
	t_r \quad 0 \quad \cdots
$$
where $t_1 \leq t_2 \leq \cdots \leq t_r \leq 0$, and the term $t_r$ occurs
in degree $r_o(C_0) + 2$.
If not all the $t_i$ are equal, let $s$, $1 \leq s \leq r$, be
the first integer for which $t_{s-1} < t_s$.  Otherwise, let $s=r+1$.
Note that $I_{C_0}$ must have elements in degree $\leq r_o(C_0)-r+3$, so we can
form the basic double link
$C_0:(1, r_o(C_0)-r+s+1) \rightarrow C_1$.
Then by Lemma~\ref{hilbfcn-bdl}, $C_1$ has Hilbert function
$$
\Delta^2 H(C_1, t) = \quad \cdots \quad t_1 + 1 \quad \cdots \quad t_{s-1}+1
	\quad t_s \quad \cdots \quad t_r \quad 0 \quad \cdots,
$$
where now the term $t_r$ occurs in degree $r_o(C_0)+3$.
Continuing by induction,
we construct a sequence of basic double links to a curve $C_m$ having
Hilbert function
$$
\Delta^2 H(C_m, t) = \quad \cdots \quad u_1 \quad \cdots \quad u_r \quad 0
\quad \cdots,
$$
where $u_1 = \cdots = u_r = 0$, and where the term $u_r$
occurs in degree $r_o(C_0) + 2 +m$.

We claim that $C_m$ has the cohomology property.  To see this, note that
$r_o(C_m) = r_o(C_0) + m$, and that
$\sigma(C_m) \leq r_o(C_0)+m-r+3$, so that $\Delta^2 H(C_m, t) = 0$
for all $t \geq r_o(C_0)+m-r+3$.  But by Corollary~\ref{Delta-coh}, we see that
$h^2(\pthree, {\cal I}_{C_m}(t)) = h^1(\pthree, {\cal I}_{C_m}(t))$
for all $t \geq r_0(C_0)+m-r-1 = r_o(C_m) - r + 1$.  That is,
$h^2(\pthree, {\cal I}_{C_m}(t)) = h^1(\pthree, {\cal I}_{C_m}(t))$
in the last $r$ places.

For the converse, suppose $C \in {\cal L}$ has the cohomology property
and that there
is a sequence of basic double links
$$
{\cal C} : C_0:(1, d_0) \rightarrow C_1:(1, d_1) \rightarrow
	\dots \rightarrow C_m = C.
$$
First,
note that by \cite[Lemma 1.14]{BM1}, $e(Z)$ increases by at least one
each time we move up in the liaison class.  But $r_o(Z)$
increases by exactly one each time.  Since the cohomology
property for $C_m$ in particular means that $e(C_m) = r_o(C_m)$,
then we must have $e(C_0) \leq r_o(C_0)$.

Next we show that the Hilbert function of $C_0$ has the
given form.  First,
$$
\Delta^2 H(C_0, r_o(C_0)+2)= \Delta^2 H(C_m, r_o(C_0) + 2 +m) -
			\delta({\cal C}, r_o(C_0)+m, m),
$$
and since $r_o(C_0) + 2 + m = r_o(C_m)+2$, we have
$\Delta^2 H(C_m, r_o(C_0)+2+m) = 0$, because of Corollary~\ref{Delta-coh}
and the
fact that $h^1(\pthree, \IC(t)) = h^2(\pthree, \IC(t))$ for $t\geq r_0(C_m)$.
Thus, $\Delta^2 H(C_0, r_o(C_0)+2) \leq 0$.  Moreover, for each
$i = 1, \dots, r-1$, we have
\begin{eqnarray*}
\Delta^2 H(C_0, r_o(C_0)+2-i) &=& \Delta^2 H(C_m, r_o(C_0)+2-i+m) -
				\delta({\cal C}, r_o(C_0)+2-i+m, m) \\
		       &=& -\delta({\cal C}, r_o(C_0)+2-i+m, m).
\end{eqnarray*}
As above, the second equality follows from the fact that
$r_o(C_0) + 2 -i + m = r_o(C_m) + 2 - i$ and because the
cohomology property for $C_m$ implies $\Delta^2 H(C_m, r_o(C_m) + 2 - i) = 0$
via Corollary~\ref{Delta-coh}.  But by our observation following
Definition~\ref{deltadef}, we have
\begin{eqnarray*}
\Delta^2 H(C_0, r_o(C_0)+2-i) &=& - \delta({\cal C}, r_o(C_0)+2-i+m, m) \\
	&\leq& -\delta({\cal C}, r_o(C_0)+2-(i-1)+m, m) \\
	&=&    \Delta^2 H(C_0, r_o(C_0)+2-(i-1)),
\end{eqnarray*}
and so the proof is finished.
\end{proof}

Since for any curve $C$ in an even liaison class ${\cal L}$ there
is a sequence of basic double links from a minimal curve
in ${\cal L}$ to a curve $C_m$, followed by a deformation
with constant cohomology to $C$, the Theorem has the immediate
consequence:

\begin{cor} \label{mincoh-prop}
An even liaison class of curves in $\pthree$ contains a curve $C$
having $h^1(\pthree, \IC(t)) = h^2(\pthree, \IC(t))$ in
the last $r$ places if and only if the minimal curve $C_0$
in ${\cal L}$ has $e(C_0) \leq r_o(C_0)$ and
has Hilbert function satisfying
$$
\Delta^2 H(C_0, r_o(C_0) - r + 3) \leq \Delta^2 H(C_0, r_o(C_0) - r + 4)
	\leq \dots \leq \Delta^2 H(C_0, r_o(C_0) + 2) \leq 0.
$$
\end{cor}

\section{Liaison Classes of Curves with Equal Cohomology}

As we saw in the previous section,
the presence or non-presence of a curve in an even liaison
class ${\cal L}$ having equal cohomology can be detected by looking at
the Hilbert function of the minimal curve in ${\cal L}$.  This raises some
immediate questions:  if an even liaison class ${\cal L}$ contains
curves with equal cohomology, in what shifts of ${\cal L}$ do they occur,
and ``how many'' such curves are there?  To answer these questions,
we first make some remarks concerning how the equal cohomology property
behaves with respect to basic double linkage.

\begin{remark}\label{eqcoh&bdls}
\begin{enumerate}
\item[(a.)] \label{bdl-rem1}
Let ${\cal L}$ be an even liaison class with associated
deficiency module of diameter $r$.
Suppose ${\cal L}$ contains curves with equal cohomology.
Then by the previous corollary, the minimal curve in the liaison class
has Hilbert function ending in a non-decreasing
sequence of $r$ non-positive terms,
beginning in degree $r_a(C_0) + 2$.
Say $\Delta^2H(C_0, t) = \cdots\quad t_1 \quad \cdots\quad t_r$ is this
sequence.  Let
$$
C_0:(1,b_0) \rightarrow C_1:(1, b_1) \rightarrow C_m:(1,b_{m-1})
	\rightarrow C_m
$$
be a sequence of basic double links.  By Lemma~\ref{hilbfcn-bdl},
if $\Delta^2 H(C_i, t)$ ends in negative terms, and if the
basic double linkage $C_i:(1,b_i) \rightarrow C_{i+1}$
changes one of these negatives, then it also changes every term
preceding.  In particular, it must change the left-most negative term.
Moreover, each basic double link which changes negatives increases
the left-most negative term by exactly one.  Note that a
link $C_i:(1, b_i) \rightarrow C_{i+1}$ changes negative terms
if and only if $b_i \geq r_a(C_i) + 2 = r_a(C_0) + i + 2$.
\item[(b.)] \label{bdl-rem2}
Related to this is the observation that if
there are more than $-t_1$ basic double links which change negative
terms, then $C_m$ cannot have equal cohomology.  More precisely,
if we have $b_i \geq r_a(C_0) + i + 2$ for more than
$-t_1$ indices $i$, then the $t_1$ term eventually becomes positive,
and this forces $h^1(\pthree, \IC(t))$ and $h^2(\pthree, \IC(t))$
to be non-equal in at least the leftmost degree.
\item[(c.)] \label{bdl-rem3}
Continuing with that theme, it follows rather trivially
that if we have more than $-t_1$ links which change
negative terms, then no further basic double link can possibly
produce a curve with equal cohomology.
\item[(d.)]  \label{bdl-rem4}
As a final remark, note that if $C :(1, d) \rightarrow D$
is a basic double link, and if $C$ has equal cohomology, then
$D$ has equal cohomology if and only if $d \leq r_a(C) + 3$. This follows
directly from Corollary~\ref{Delta-coh} and Lemma~\ref{hilbfcn-bdl}.
\end{enumerate}
\end{remark}

Now, we begin our description of which curves in the liaison
class have equal cohomology by showing that there is a
unique minimal such curve.

\begin{prop} Suppose ${\cal L}$ is an even liaison class which
contains curves having equal cohomology.  Then up to deformation
through curves with constant cohomology, there is a unique
minimal curve with equal cohomology.
\end{prop}

\begin{proof}  This is quite easy, given our remarks above.
If $\Delta^2 H(C_0, t) = \cdots\quad t_1\quad\cdots\quad t_r$ are
the final non-decreasing non-positive terms in the Hilbert function of $C_0$,
as guaranteed by Corollary~\ref{mincoh-prop},
then it requires at least $-t_1$ basic double links to reach
a curve $D$ for  which $\Delta^2 H(D, t) = 0$ for $t \geq r_a(D) + 2$,
and by Corollary~\ref{Delta-coh}, $D$ then has equal cohomology.
On the other hand, the construction in Theorem~\ref{coh-prop}
produces a curve with equal cohomology in exactly $-t_1$ steps.
Also, if we reach a curve in equal cohomology in $-t_1$ steps,
then every basic double link of degree $b_i$ satisfies
$b_i \geq r_a(C_0) + i + 2$.  Thus, by Lemma~\ref{hilbfcn-bdl},
no matter what sequence of basic double links we take, the Hilbert
function of the resulting curve is invariant.  Thus
any two curves in this shift with equal cohomology are deformations
through curves with constant cohomology;  see \cite[Proposition 3.1]{BM2}.
\end{proof}

Now, we need to recall a result from \cite{BM2} on ``equivalence''
of basic double links.  In the proof of \cite[Lemma 5.2]{BM2},
they show the following:  suppose
$$
C_1:(1, b_1) \rightarrow C_2:(1, b_2) \rightarrow C_3
$$
are basic double links with $b_1 < b_2$.
Then the sequence $b_1, b_2$ is equivalent
to the sequence $b_2 - 1, b_1 + 1$ in the sense that if we make
the basic double links
$$
C_1:(1, b_2 - 1) \rightarrow C_2':(1, b_1+1) \rightarrow C_3',
$$
then $C_3$ and $C_3'$ are deformations of each other, through
curves with constant cohomology.  Note that implicit in this
is the fact that the basic double linkage of degree $b_2 - 1$
can actually be made; this is also noted in their proof.

We will use this idea of ``flipping'' adjacent degrees in the
next result.

\begin{prop}  Suppose ${\cal L}$ is an even liaison
class containing curves with equal cohomology, and that
$C$ is a curve in ${\cal L}$ having equal cohomology.
Then $C$ is obtained from the minimal curve having
equal cohomology by a sequence of basic double linkages
followed by a deformation through curves with constant cohomology,
if necessary. Each curve in the sequence also has equal cohomology.
\end{prop}

\begin{proof}  Assume that the deficiency module associated
to ${\cal L}$ has diameter $r$.  First, since $C$ is in ${\cal L}$, then
by the Lazarsfeld--Rao property, after deforming $C$ through curves
with constant cohomology if necessary, we may assume that
there is a sequence of basic double links
\begin{equation}\label{bdl1}
C_0:(1,b_0) \rightarrow C_1:(1,b_1) \rightarrow\cdots \rightarrow C_m = C
\end{equation}
where $C_0$ is the (absolute) minimal curve in ${\cal L}$,
and we can assume that $b_0 \leq \dots \leq b_{m-1}$.
Since ${\cal L}$ possesses curves with equal cohomology,
$C_0$ has Hilbert function given by
$$
\Delta^2 H(C_0, t) = \quad\cdots\quad t_1 \quad\cdots\quad t_r
$$
where $t_1\leq\dots\leq t_r \leq 0$, and $t_1$ is in
degree $r_a(C_0) + 2$.
Now, since $C$ has equal cohomology, then exactly $-t_1$
of the basic double links in the sequence~(\ref{bdl1})
change negatives.  That is, exactly $-t_1$ of the $b_i$
satisfy $b_i \geq r_a(C_0) + i + 2$.  This follows from
Remark~\ref{bdl-rem2}(a).

Choose the first index $s \geq 0$ for which $b_i \geq r_a(C_0) + i + 2$,
and using the equivalence outlined above, flip this degree down
to the first position.  Note that
this is possible since our original sequence of $b_i$'s is non-decreasing
and $b_{s-j} < r_a(C_0) + s + 2 - j \leq b_s - j$
for $0 \leq j \leq s$.
This creates an equivalent sequence
of basic double links of degrees $b_0', \dots, b_{m-1}'$,
where $b_0' = b_s - s$, $b_i' = b_{i-1} + 1$ for $0 < i \leq s$,
and $b_i' = b_i$ for $i \geq s+1$.  In particular, exactly
$-t_1$ of the $b_i'$ satisfy $b_i' \geq r_a(C_0) + i + 2$,
and moreover $b_0' \geq r_a(C_0) + 2$.  Continue in the same manner:
find the second time that a $b_i'$ changes negatives, and flip
it down to the second position, and so forth.  Then we end up
with a sequence $c_0, \dots, c_{m-1}$ of basic double links
which is equivalent to the one we started with, and which
moreover has $c_i \geq r_a(C_0) + i + 2$ for $i = 0, \dots, -t_1-1$.
Hence, by Remark~\ref{eqcoh&bdls}, since we change exactly
$-t_1$ negative terms, all in the first $-t_1$ links, and since
we eventually end up with a curve having equal cohomology,
then the curve $C_{-t_1}$ and each curve from $C_{-t_1}$ on,
must also have equal cohomology.  This follows from Remark~\ref{bdl-rem1}.
\end{proof}

We recapitulate what we have proven in the next statement:

\begin{thm} \label{LR-prop}
Suppose ${\cal L}$ is an even liaison class
containing curves with equal cohomology.  Then there is a minimal
shift ${\cal L}^t$ which contains a curve with equal cohomology;
the curves with equal cohomology in the minimal shift are
unique up to deformation
through curves with constant cohomology; every curve in ${\cal L}$
with equal cohomology is obtained from the minimal one
by basic double linkage and deformation through curves with
constant cohomology; and finally
every rightward shift of ${\cal L}^t$ contains, up
to deformation, a finite, non-zero number of curves
with equal cohomology.
\end{thm}

\begin{proof} The only part which remains to be proven is
the final statement.  If $C$ is a curve with equal cohomology
in some shift ${\cal L}^s$ of the liaison class, then
$\Delta^2 H(C, t) = 0$ for all $t \geq r_a(C) + 2$.
This implies in particular that $I_C$ contains non-zero
elements of degree $\geq r_a(C) + 2$.  Thus, we can make a
basic double link $C:(1, r_a(C) + 2) \rightarrow D$,
and $D \in {\cal L}^{s+1}$ also has equal cohomology.
On the other hand,  if we make a basic double link
of degree $> r_a(C) + 3$, then the resulting curve does not
have equal cohomology, so there is only a finite number
of allowable degrees.
\end{proof}

Theorem~\ref{LR-prop} shows that the curves with equal
cohomology have a strong Lazarsfeld-Rao property,
in the sense that there are unique minimal curves,
every other curve is obtained by basic double linkage, and
in each allowable shift, there are only a finite number
of curves, up to deformation.

In the case of Buchsbaum liaison classes, we can actually count
the number of curves in each shift which have equal cohomology.

\begin{prop} \label{number-Buchsbaum}  Suppose ${\cal L}$ is a
Buchsbaum even liaison class having curves with equal cohomology,
and let ${\cal L}^s$ be the minimal shift in which such
a curve occurs.  Then for each $t \geq s$, there are
exactly $2^{t-s}$ curves, up to flat deformations, having equal cohomology.
\end{prop}

\begin{proof}  This follows from the fact that if $D \in {\cal L}^h$,
then $\alpha(I_D) = 2N + h$, where $N = \sum \dim H^1_*(\pthree, {\cal
I}_D(i))$,
and the description of the Hilbert function of minimal Buchsbaum
curves in \cite[Proposition 2.1]{BM1}.
In particular, the minimal curve $C$ having equal cohomology has
Hilbert function
$$
\Delta^2 H(C, t) = 1 \quad 2 \quad \cdots \quad 2\alpha + s.
$$
In order to move from ${\cal L}^s$ to ${\cal L}^{s+1}$ and preserve
the cohomology property, we can only make basic double
links of degree $2\alpha+s$ or $2\alpha+s+1$.
Similarly, we can only move from ${\cal L}^{s+1}$ to ${\cal L}^{s+2}$
by basic double links of degree $2\alpha+s+1$ or $2\alpha+s+2$.
Continuing inductively, the statement is proven.
\end{proof}

We can also give a proof more along the lines of the original
proof that the liaison classes of curves in $\pthree$ have
the Lazarsfeld--Rao property.  It is much less constructive
in nature, but, in some sense, points out the naturality
of our cohomological criterion of equal cohomology.
Since it is so different in spirit from our previous argument, we
felt it necessary to include it here.

\begin{lemma} Let $\cal F$ be a rank $(r+1)$ vector bundle on $\pthree$ with
$H^2_* (\pthree, {\cal F}) = 0$.  Let
\[
\phi_1 : \bigoplus_{i=1}^{r+1} {\cal O}_{\pthree} (-a_i ) \rightarrow
{\cal F},
\hskip .5in a_1 \leq \dots \leq a_r
\]
\[
\phi_1 : \bigoplus_{i=1}^{r+1} {\cal O}_{\pthree} (-b_i ) \rightarrow
{\cal F},
\hskip .5in b_1 \leq \dots \leq b_r
\]
be morphisms whose degeneracy loci are curves $C_1$ and $C_2$ with equal
cohomology.  Then there exists a morphism
\[
\phi : \bigoplus_{i=1}^{r+1} {\cal O}_{\pthree} (-c_i ) \rightarrow {\cal F},
\hskip .5in c_i = \min \{ a_i ,b_i \}
\]
whose degeneracy locus is also a curve with equal cohomology.

\end{lemma}

\begin{proof}
Notice that $C_1$ and $C_2$ are evenly linked, in the even liaison class
determined by the stable equivalence class of $\cal F$, according to Rao's
classification \cite{rao}.  By \cite[Lemma 2.1]{BBM}, there exists such a
$\phi$ whose degeneracy locus is a curve
$C$.  We just have to prove that $C$ has equal cohomology.

Twisting and relabeling if necessary, we may assume that we have locally free
resolutions
\[
0 \rightarrow \bigoplus_{i=1}^{r+1} {\cal O}_{\pthree} (-a_i ) \rightarrow
{\cal
F} \rightarrow {\cal I}_{C_1} \rightarrow 0
\]
\[
0 \rightarrow \bigoplus_{i=1}^{r+1} {\cal O}_{\pthree} (-b_i ) \rightarrow
{\cal
F} \rightarrow {\cal I}_{C_2} (h) \rightarrow 0.
\]
Notice that the deficiency module of $C_2$ is shifted $h$ places to the right
of that of $C_1$.  We then get
\[
0 \rightarrow H^2 (\pthree, {\cal I}_{C_1} (t)) \rightarrow
	H^3 (\pthree, \bigoplus {\cal O}_{\pthree} (t-a_i )) \rightarrow
	H^3 (\pthree, {\cal F}(t)) \rightarrow 0
\]
\[
0 \rightarrow H^2 (\pthree, {\cal I}_{C_2} (t+h)) \rightarrow
	H^3 (\pthree, \bigoplus {\cal O}_{\pthree} (t-b_i )) \rightarrow
	H^3 (\pthree, {\cal F}(t)) \rightarrow 0.
\]
(The first 0 comes from the assumption on the vanishing of the cohomology of
$\cal F$ and the second from the fact that $h \geq 0$.)  By the assumption of
equal cohomology, for $t \geq r_a(C_1)$
the first term in the first sequence has the
same dimension as the first term in the second sequence.
Hence for $t \geq r_a(C_1)$, also the second terms are equal.  Therefore
\[
\left \{ a_i \ | \ r_a(C_1) - a_i \leq -4 \right \} =
\left \{ b_i \ | \ r_a(C_1) - b_i \leq -4 \right \}
\]
(since these are the terms which contribute to the middle cohomology space in
the degrees $t \geq r_a$).  That is,
\[
\left \{ a_i \ | \ a_i \geq r_a + 4 \right \} =
\left \{ b_i \ | \ b_i \geq r_a + 4 \right \}
\]
Call this set $A$.

Now, for any curve $Y$ with locally free resolution
\[
0 \rightarrow \bigoplus_{i=1}^{r+1} {\cal O}_{\pthree} (-d_i ) \rightarrow
{\cal F} \rightarrow {\cal I}_{Y} (\delta ) \rightarrow 0,
\]
$Y$ has equal cohomology if and only if
$\{ d_i \ | \ d_i \geq r_a(Y) + 4 \} = A$
(since this set determines $h^3(\pthree, \bigoplus {\cal O}(-d_i + t))$
and hence $h^2 ({\cal I}_Y (t+\delta )$ in the desired range).
The proof of the lemma follows immediately from this fact.
\end{proof}

\begin{cor} \label{LRprop}
The set of curves in a given even liaison class which have equal cohomology
satisfy the Lazarsfeld--Rao property.
\end{cor}

\begin{proof}
The proof is identical to that in \cite{BBM}.  The lemma above replaces
\cite[Lemma 2.1]{BBM}.    Then \cite[Proposition 2.2]{BBM}
goes through to prove
the uniqueness of the minimal element.  Similarly, \cite[Proposition 2.3]{BBM}
goes through to show the relation between the minimal element and any other
curve in the even liaison class with equal cohomology.  Finally,
\cite[Theorem 2.4]{BBM}
still works to show how to produce a curve with equal cohomology as
a sequence of basic double links followed by a deformation, starting with a
minimal curve with equal cohomology.  The proof in \cite{BBM} shows that such a
sequence exists.  The fact that we start with equal cohomology and end
with equal cohomology shows that every step in between has equal cohomology
too.
\end{proof}

\begin{remark}\label{mincurves}  We do not yet know of any examples
of even liaison classes for which the absolute minimal curve
$C_0$ is also the minimal curve with equal cohomology.
\end{remark}

\section{Integral Curves with Equal Cohomology} \label{integral}
There has been much recent progress on further clarifying the
structure of even liaison classes by giving conditions for the
presence within the liaison classes of nice curves.  In particular,
there is information on where in a given class one
can find integral curves \cite{nollet},
or smooth and connected curves in Buchsbaum classes \cite{MDP2}, and
on how these curves are related to each other and to the minimal
curve in the class.
The paper \cite{PR} shows that at least in Buchsbaum
classes, smooth and connected curves share the same
Lazarsfeld--Rao properties as irreducible curves, and
imply that one can obtain the integral curves within
a given shift of a liaison class by deforming irreducible curves.
Their calculations are based also on the work of Nollet, as well
as on \cite{MDP2}.

In this section, we are interested in using some of the results
of \cite{nollet} to obtain some information on
when an even liaison class contains curves which are integral
and have equal cohomology.  We first recall the relevant
definitions and some results from \cite{nollet}.

Let $C$ be a curve in $\pthree$, defined by an ideal $I=I_C$.
The postulation character of $C$ is given by $\gamma_C(n) = -\Delta^3 H(C, n)$.
There are three natural invariants to attach to $C$:
\begin{eqnarray*}
s(C) &=& \min\{\,n : \gamma_C(n) \geq 0 \,\} \\
t(C) &=& \min\{\,n : \gamma_C(n) > 0\,\} \\
t_1(C) &=& \mbox{smallest degree of a surface containing $C$ which meets}\\
&&\quad\quad\mbox{a surface of degree $s(C)$ containing $C$ properly}.
\end{eqnarray*}

We note for clarity that $s(C) = \alpha(C)$, the minimal degree of
a generator of $I_C$, and $t_1(C) = \beta(C)$, the minimal degree
for which $I_{\leq t} = \oplus_{i \le t} [I]_i$ generates
an ideal of codimension $2$.
Next, say that $C$ dominates a curve $D$ at height $h$ if
$C$ can be obtained from $D$ by a sequence of $h$ basic double
links, followed by a deformation.  The central definition for this
section is the following:
suppose $C$ dominates the minimal curve $C_0$ in ${\cal L}$ at height $h$.
Then
$$
\theta_C(n) = \left\{ \begin{array}{ll}
			\gamma_C(n), & \mbox{if $s(C) \le n < s(C_0) + h$} \\
			\gamma_C(n) - \gamma_{C_0}(n-h), &
				      \mbox{if $n \geq s(C_0) + h$} \\
			0,	    & \mbox{otherwise.}
		     \end{array} \right.
$$
(This definition appears in \cite{PR} and is clearly equivalent
to the one in \cite{nollet}.)
We say $\theta_C$ is connected in degrees $\ge a$ if
$\theta_C(b) > 0$ for $b \geq a$ implies $\theta_C(n) > 0$
for all $a \leq n \leq b$, and similarly $\theta_C$ is connected
in degrees $\le b$ if $\theta_C(a) > 0 $ for some $a \leq b$
implies $\theta_C(n) > 0$ for all $a \leq n \leq b$.  Finally,
$\theta_C$ is connected about an interval $[a,b]$ if it
is connected in degrees $\geq a$ and in degrees $\leq b$,
and if $\theta_C(n) > 0$ for all $a \leq n \leq b$.

Now, Nollet proves the following theorem in \cite{nollet}:

\begin{thm}  Let ${\cal L}$ be an even liaison class of curves
in $\pthree$ with minimal curve $C_0$.
\begin{enumerate}
\item[{\rm (}a.{\rm )}] {\rm (\cite[Theorem 5.2.1]{nollet})}
If $C \in {\cal L}$ is an integral curve of height $h$,
then $\theta_C$ is connected about $[t(C_0)+h, t_1(C_0) + h - 1]$.
\item[{\rm (}b.{\rm )}] {\rm (\cite[Theorem 5.2.5]{nollet})}
Conversely, suppose
$C$ dominates at height $h$ an integral curve $D$ in ${\cal L}$,
which is generically Cartier on a surface of minimal degree and
has either $\theta_D \not=0$ or $t(D) \leq e(D) + 4$.  If
$\theta_C$ is connected about $[t(C_0) + h, t_1(C_0) + h - 1]$, then
$C$ can be deformed to an integral curve.
\end{enumerate}
\end{thm}

Thus, having $\theta_C$ connected about the interval
$[t(C_0) + h, t_1(C_0) + h - 1]$ is very close to having $C$ integral.
As it turns out, this condition is relatively easy to check
for curves with equal cohomology.  We begin with some
elementary calculations.  Throughout, let ${\cal L}$ be an
even liaison class of curves, which contains curves having equal
cohomology, and let $C_0$ be the (absolute) minimal curve in ${\cal L}$
and $C$ the minimal curve in ${\cal L}$ with equal cohomology.

\begin{lemma}\label{s_t_and_theta}
Suppose the minimal curve $C$ with equal cohomology has height $h$
over $C_0$.  Then:
\begin{eqnarray*}
s(C) &=& s(C_0) + h \\
t(C) &=& t(C_0) + h \\
\theta_C(n) &=& \left\{ \begin{array}{ll}
   -\gamma_{C_0}(n-h) & \mbox{\rm if $r_a(C_0) + h + 2 < n \leq \sigma(C_0) +
h$}\\
    0 & \mbox{\rm otherwise }.\end{array}\right.
\end{eqnarray*}
\end{lemma}

\begin{proof}  Note that in order to move up the liaison class from the
minimal curve $C_0$ to the curve $C$, in order to get $C$ with equal
cohomology, we must take basic double links $C_i \rightarrow C_{i+1}$
of degree large enough to change the negative signs in $\Delta^2 H(C_i, t)$.
Clearly, this degree $d$, say, is strictly larger than $t(C_i)$.
By \cite[Corollary 2.3.5]{nollet}, then,
$\gamma_{C_{i+1}}(n) = \gamma_{C_i}(n-1)$ for $n \leq d$.  In particular,
$s(C_{i+1}) = s(C_i) + 1$ and $t(C_{i+1}) = t(C_i) + 1$, and so the
first two statements are done by induction.  The assertion about
$\theta_C$ then follows from the definition of $\theta_C$ and
the fact that $\gamma_C(n) = 0$ for $n \geq r_a(C_0) + h + 2$,
since $C$ has equal cohomology, and using that $s(C_0) \leq r_a(C_0) + 2$
and $s(C) = s(C_0) + h$.
\end{proof}

Now we are ready to determine which curves have both equal cohomology
and connected $\theta$.  Our first result takes care of a rather
trivial case.

\begin{prop}  Suppose the minimal curve $C$ with equal cohomology has
height $h$ over $C_0$ and has $\theta_C$ connected about the interval
$[t(C_0) + h, t_1(C_0) + h - 1]$.  Then $C=C_0$, up to deformation,
and $t(C_0) = t_1(C_0)$.

Conversely, if $C_0$ has equal cohomology and $t(C_0) = t_1(C_0)$,
then $\theta_{C_0}$ is connected about the interval $[t(C_0), t_1(C_0) - 1]$.
\end{prop}

\begin{proof}  First note that since $C$ has equal cohomology,
we must have $t(C) \leq r_a(C) + 2$, and this then implies by
Lemma~\ref{s_t_and_theta} that $t(C_0) \leq r_a(C_0) + 2$.
But again by Lemma~\ref{s_t_and_theta}, this means that
$\theta_C(t(C_0)+h) = 0$, so $\theta_C > 0$ on the
interval $[t(C_0)+h, t_1(C_0) + h -1]$ if and only if this
interval is empty.  Clearly, this is equivalent to $t(C_0) = t_1(C_0)$.
Now, $\theta_C$ is connected in degrees $\geq t(C_0) + h$ if and
only $\theta_C = 0$ in degrees $\geq t(C_0) + h$, and again
by Lemma~\ref{s_t_and_theta}, this occurs if and only if
$\sigma(C_0) = r_a(C_0) + 2$.  By Corollary~\ref{Delta-coh},
this implies that $C_0$ has equal cohomology, and so $C = C_0$,
up to deformation.

The other direction is quite trivial, since $\theta_{C_0} = 0$ and
$[t(C_0), t_1(C_0) - 1]$ is empty.
\end{proof}

Our next proposition is the main result of this section, and tells
us when an even liaison class contains curves with
equal cohomology and connected $\theta$.  Note that it is identical
in spirit to Theorem~\ref{mincoh-prop}.

\begin{prop}  \label{theta-coh}
Suppose ${\cal L}$ is an even liaison class of curves with minimal
curves satisfying $t(C_0) < t_1(C_0)$.  Then ${\cal L}$ contains
a curve $D$ of height $\hover$, say, having equal cohomology
and $\theta_D$ connected about $[t(C_0) + \hover, t_1(C_0) + \hover - 1]$
if and only if $t_1(C_0) \leq \sigma(C_0) + 1$ and the Hilbert function
of $C_0$ satisfies
$$
\Delta^2H(C_0, r_a(C_0) + 2) < \Delta^2H(C_0, r_a(C_0) + 3) < \cdots
	<\Delta^2H(C_0, \sigma(C_0) - 1) < 0.
$$
\end{prop}

\begin{proof}  First, suppose $C_0$ has $t_1(C_0) \leq \sigma(C_0) + 1$
and satisfies the condition on the Hilbert function.   Note
that the condition on the Hilbert function implies
that $t(C_0) \leq r_a(C_0) + 2$.  Then by
Lemma~\ref{s_t_and_theta}, the minimal curve $C$ with equal
cohomology has $\theta_C(t) = 0$ for $t(C_0) + h \leq t \leq r_a(C_0) + h + 2$
and $\theta_C(t) > 0$ for $r_a(C_0)+h+2 < t \leq \sigma(C_0) + h$,
where $C$ has height $h$ over $C_0$.
Now perform a sequence of $r_a(C_0) + 2 - t(C_0)$ basic double links
all of degree $r_a(C_0) + h + 3$ to reach a curve $D$ of
height $\hover = h + r_a(C_0) + 2 - t(C_0)$.  Then by a repeated
application of \cite[Corollary 2.3.5]{nollet}, it is easy to
check that
$$
\theta_D(t) = \left\{
\begin{array}{ll}
1 & \mbox{ for $t(C_0) + \hover \leq t \leq r_a(C_0) + \hover + 2$} \\
\theta_C(t+h-\hover) & \mbox{ for $r_a(C_0) + \hover + 2 < t \leq \sigma(C_0) +
\hover$} \\
0 & \mbox{ otherwise. }
\end{array} \right.
$$
In particular, $\theta_D$ is connected about
$[t(C_0) + \hover, t_1(C_0) + \hover - 1]$.  Also, since $D$ was
obtained from $C$ by basic double links of low degree, $D$ still
has equal cohomology; see Remark~\ref{bdl-rem4}(d.).

Conversely, suppose $D$ is a height $\hover$ curve with equal cohomology,
and with $\theta_D$ connected about the interval
$[t(C_0) + \hover, t_1(C_0) + \hover - 1]$.  Then there is a sequence
of basic double links
$$
C:(1,b_0) \rightarrow C_1:(1, b_1) \rightarrow \cdots
	\rightarrow D
$$
where $C$ is the minimal curve in ${\cal L}$ with equal cohomology
(of height $h$, say) and where each $b_i$ satisfies
$b_i \le r_a(C_0) + h + i + 3$ (see Remark~\ref{eqcoh&bdls}(d.) or
the proof of Corollary~\ref{LRprop}).  If $C=C_0$, then
$r_a(C_0) + 2 = \sigma(C_0)$, so there is nothing
to show.  Hence we may assume that $C_0$ does not have
equal cohomology, and this means that $\Delta^2 H(C_0, r_a(C_0) + 2) < 0$,
which in turn implies $\gamma_{C_0}(r_a(C_0) + 3) < 0$
and $t(C_0) \leq r_a(C_0) + 2$.

By a repeated use of \cite[Corollary 2.3.5]{nollet},
$\gamma_D(n) = \gamma_C(n-\hover+h)$ for $r_a(C_0) + \hover+3 \leq n$.
But by Lemma~\ref{s_t_and_theta},
$\gamma_C(n-\hover+h) = -\gamma_{C_0}(n-\hover)$ for
$r_a(C_0) + \hover + 3 \leq n \leq \sigma(C_0) + \hover$.

Now, $\theta_D$ is connected in degree $\geq t(C_0) + \hover$, and
our assumption that $t(C_0) < t_1(C_0)$ implies in particular that
$\theta_D(t(C_0) + \hover) > 0$.
Also, $\theta_D(\sigma(C_0) + \hover) = -\gamma_{C_0}(\sigma(C_0)) > 0$,
so $\theta_D > 0$ on the interval $[t(C_0) + \hover, \sigma(C_0) + \hover]$.

Next, note that $\theta_D$ is positive on the
interval $[r_a(C_0) + \hover + 3, \sigma(C_0) + \hover]$, since this
interval is contained in the interval $[t(C_0) + \hover, \sigma(C_0) +
\hover]$.
Thus, we have $0 < \theta_D(t) = -\gamma_{C_0}(t-\hover)$
for $r_a(C_0) + \hover + 3 \leq t \leq \sigma(C_0) + \hover$.
This clearly implies that $\gamma_{C_0}(t) < 0$ for
$r_a(C_0) + 3 \leq t\leq \sigma(C_0)$, and this means
that $\Delta^2 H(C_0, t)$ is strictly increasing in the given range.

Finally, to see that $t_1(C_0) \leq \sigma(C_0) + 1$, note
that the connectedness property of $\theta_D$ implies that
$\theta_D > 0$ on $[r_a(C_0) + \hover + 3, t_1(C_0) + \hover]$,
but clearly $\theta_D(t) = 0$ for $t >\sigma(C_0) + \hover$,
by the argument above.
\end{proof}

\section{Degrees of Generators and Liaison Classes of Curves}

In this section, we go back to studying degrees of generators of
the ideals defining space curves, by using the results on cohomology
given in the previous sections.  We are able to give some nice conditions
on the degrees of the components of the deficiency module associated
to a liaison class in order for the class to contain a curve whose
ideal has a generator of maximal degree.  Since knowledge about curves
with equal cohomology in a given liaison class depends so crucially
on knowing the Hilbert function of the minimal curve in the liaison
class, our characterizations for when a liaison class contains
curves with generators of high degree work best when we already
know the Hilbert function of the minimal curve.  To do this, we
have to make some extra assumptions on the liaison class.
We have concentrated on cohomological criteria, and two
very clean statements are given in Proposition~\ref{degree-maxcorank}
and in Proposition~\ref{degree-Buchsbaum}.  Similarly, using
our results on integral curves, we give some results on
existence within a liaison class of integral curves with
generators of maximal degree.

As we showed in the previous sections, if $C$ is defined by
an ideal $I = I_C$, then the property of $I_C$ having
a generator of high degree is very closely related to
having $h^1(\pthree, \IC(t)) = h^2(\pthree, \IC(t))$ in
a large number of places, and furthermore, curves with
this cohomology property are easily constructed by basic double linkage,
as long as we know the minimal curve in a liaison class.  However,
we should remark that in order to construct a curve whose
ideal has a high degree generator, we need to choose the degrees
of the basic double links carefully, since it is possible to
make $h^1 = h^2$ in the maximum number of places, without
introducing a high degree generator.  The following example
should clarify this somewhat.

\begin{example}  Start with the Buchsbaum liaison class
${\cal L} = {\cal L}_{4,1}$, whose deficiency module has two
consecutive components of dimensions $4$ and $1$, respectively.
Then the minimal curve $C_0$ in ${\cal L}$ has
Hilbert function (see \cite[Corollary 2.18]{BM1})
$$
\Delta^2 H(C_0, t) = 1 \quad 2 \quad \cdots \quad 10 \quad -2 \quad -1.
$$
Taking the two basic double links
$$
C_0:(1, 13) \rightarrow C_1:(1, 13) \rightarrow C_2
$$
produces the curve $C_2$, for which $h^1 = h^2$ in the last two places,
and which has $\sigma(C_2) = 12$, but
$\omega(I_{C_2}) = 13$. Thus the bound in Proposition~\ref{maxdegree}
is not obtained, even though $C_2$ does have the cohomology property.

On the other hand, note that if we take the basic double links
$$
C_0:(1, 12) \rightarrow C'_1:(1,14) \rightarrow C'_2,
$$
then $C'_2$ has $h^1 = h^2$ in the last two places, and
also has $\sigma(C'_2) = 12 = \omega(I_{C'_2}) - 2$. So the
bound is achieved in this case.
\end{example}

Generally speaking, this second sequence of basic double links
produces curves whose defining ideals have a high degree generator.
Note that it is essentially the procedure given by Theorem~\ref{coh-prop}.
However, we need to require that there are no ``trailing zeroes''
on the end of the second difference of the Hilbert function.  We
formalize this in the next statement.

\begin{prop}\label{degree}  Suppose ${\cal L}$ is an even liaison
class of curves in $\pthree$, with minimal curve $C_0$, and
let $r = \diam H^1_*(\pthree, \IC)$. If the Hilbert function
of $C_0$ satisfies
$$
\Delta^2 H(C_0, r_a(C_0)+2) \leq \dots \leq \Delta^2 H(C_0, r_o(C_0)+2) < 0,
$$
then there exists a curve $C$ in ${\cal L}$ whose
defining ideal $I_C$ satisfies $\omega(I_C) = \sigma(S/I_C) + r$.
\end{prop}

\begin{proof}
We follow the construction of basic double links in the
first part of Theorem~\ref{coh-prop}.  Perform the given
sequence of basic double links, up to the next to the last stage. Since
$\Delta^2 H(C_0, r_o(C_0) + 2) < 0$, then at this stage
we see that
$$
\Delta^2 H(C_{m-1}, t) = \cdots \quad -1 \quad -1 \quad \cdots \quad -1 \quad 0
	\quad \cdots,
$$
where there are $r$ terms of $-1$, and the rightmost term
occurs in degree $r_o(C_0)+2+m$.  Thus our final basic double link
$C_{m-1}:(1, r_o(C_0)+m+4) \rightarrow C_m$ produces
the curve $C_m$, whose ideal $I_{C_m}$ has a generator
of degree $r_o(C_0) + m + 4$, and such that
$\sigma(C_m) = r_o(C_0) + m - r + 4$.  Hence
$\omega(I_{C_m}) = \sigma(C_m) + r$, which is what we wanted
to show.
\end{proof}

In fact, using the structure theory for curves with equal cohomology
developed in the previous section, we can say more about the
curves in an even liaison class having maximal degree generators.
Namely, the minimal such curve occurs in the shift
${\cal L}^s$ of ${\cal L}$, where $s=-\Delta^2 H(C_0, r_a+2)$, and
every other such curve is obtained by a sequence of basic double
links from this minimal one, followed if necessary by a deformation.
Moreover, up to a flat deformation, there are only a finite number
of curves with a maximal degree generator in each allowable shift.

\begin{remark}  It is interesting to note that curves in a
non-arithmetically Cohen--Macaulay liaison class which
have a generator of maximal degree in fact have all
of their generators of relatively high degree.  Indeed,
since the minimal such curve $C$ lies in the shift ${\cal L}^s$
as above, $\alpha(I_C) = \alpha(I_{C_0}) + s$.
So, at least if the minimal curve does not already have
a maximal degree generator (see Remark~\ref{mincurves}),
we are forced to increase all the degrees of the generators.

This is in some contrast with the arithmetically Cohen--Macaulay
case, where there are curves $D$ with quadric generators
and for which $\omega(D) = \sigma(D)$, the maximum possible.
For example, let $D$ be the union of a plane curve of
degree $m$
and a line, which meet at one point.  Then it is easy to see
that $I_D$ has quadric generators and $\omega(D) = \sigma(D) = m$.

On the other hand, even within a non-arithmetically Cohen--Macaulay
liaison class, we can make the difference $\omega(I_C) -\alpha(I_C)$
arbitrarily large when $\omega(I_C)$ is maximal.  For this,
simply take the minimal curve $C$ with a maximal degree
generator, and form $m$ basic double links all of
degree $\alpha(I_C)$.  Then the resulting curve $C_m$
still has $\omega(I_{C_m})$ maximal (see Remark~\ref{eqcoh&bdls}(d.)),
and has
$\omega(I_{C_m}) - \alpha(I_{C_m}) = \omega(I_C) - \alpha(I_{C}) + m$.
\end{remark}

Next, we want to interpret the conditions on Hilbert functions
only in terms of cohomology.  As we noted above, to do this
we need to make some extra assumptions to allow us to calculate
the Hilbert function of the minimal curve.  Our first result
is for maximal corank curves, and the second for Buchsbaum curves.

\begin{cor} \label{degree-maxcorank}
Suppose ${\cal L}$ is an even liaison
class of curves in ${\Bbb P}^3$, with minimal curve $C_0$.  Assume that
$e(C_0) < r_a(C_0)$ {\rm (}i.e., $C_0$ has maximal corank{\rm )}.
Then ${\cal L}$ contains a curve $C$ such that
$\omega(I_C) = \sigma(C) + \diam H^1_*(\pthree, \IC)$ if and only if
$$
h^1(\pthree, {\cal I}_{C_0}(t)) \geq 3 h^1(\pthree, {\cal I}_{C_0}(t+1))
		- 3 h^1(\pthree, {\cal I}_{C_0}(t+2))
		+   h^1(\pthree, {\cal I}_{C_0}(t+3)),
$$
for $t=r_a(C_0), \dots, r_o(C_0)$.
\end{cor}

\begin{proof}  The conditions on the cohomology modules
guarantee via Lemma~\ref{hilbfcn_eqn} that
$$
\Delta^2 H(C_0, r_o - r + 2) \leq \dots \leq \Delta^2(r_o + 2, C_0)
					< 0,
$$
and so sufficiency follows from Proposition~\ref{degree}.

To see necessity, note that if ${\cal L}$ contains a curve
with a generator of maximal degree, then by Proposition~\ref{equal-coh},
that curve has equal cohomology, and so by Corollary~\ref{mincoh-prop},
the minimal curve in ${\cal L}$ has Hilbert function whose
second difference ends in a sequence of non-decreasing negative
terms, and then Lemma~\ref{hilbfcn_eqn} translates this
back to the required statement about cohomology.
\end{proof}

Recall that a Buchsbaum liaison class is completely determined
by the dimensions of the graded components of the associated deficiency module.
We will write ${\cal L}_{n_1\dots n_r}$ for the Buchsbaum class
associated to the graded module $M = \oplus [M]_i$, where
$\dim_k [M]_i = n_i$ and is zero otherwise, and where we assume that
$n_1, n_r > 0$ and $n_i \ge 0$ for $1 < i < r$.

\begin{prop} \label{degree-Buchsbaum}
Suppose ${\cal L} = {\cal L}_{n_1 \dots n_r}$ is a Buchsbaum even liaison
class.  Then ${\cal L}$ contains a curve $C$ such that
$\omega(I_C) = \sigma(I_C) + r$ if and only if
$n_i \geq 3 n_{i+1}$ for $i = 1, \dots r-1$.
\end{prop}

\begin{proof}  It follows from \cite[Corollary 2.18]{BM1} that
the conditions on the deficiency module guarantee that the
Hilbert function of $C_0$ satisfies the conditions given in
Proposition~\ref{degree}, and we can therefore use that result to
prove sufficiency.

On the other hand, if $C$ is a curve in ${\cal L}$ whose ideal has a generator
of maximal degree, then by Proposition~\ref{equal-coh}, $C$ has
equal cohomology.
Thus by Corollary~\ref{mincoh-prop}, the
minimal curve $C_0 \in {\cal L}$ has Hilbert function
satisfying
$$
\Delta^2 H(C_0, r_o(C_0) - r + 2) \leq \dots \leq \Delta^2 H(C_0, r_o(C_0) + 2)
			\leq 0.
$$
But now it follows from the description of the Hilbert function for minimal
Buchsbaum curves given in \cite[Corollary 2.18]{BM1} that
$n_i \geq 3 n_{i+1}$, for each $i = 1, \dots, r-1$.
\end{proof}

\begin{remark}
Proposition \ref{degree} is, in a sense, a complete answer to the problem of
determining which even liaison classes $\cal L$ contain a curve $C$ whose
defining ideal $I_C$ satisfies
$\omega (I_C ) = \sigma (S/I_C ) + \hbox{diam } H^1_* (\pthree, {\cal I}_C )$.
However, it is generally not easy to tell, for a given even liaison class,
what the Hilbert function of the corresponding minimal
curve is.  So it is worth noting that there are also times when one can tell
directly from the associated deficiency module $M$ (defined up to shift)
that such a curve does not exist.

If $M$ is annihilated by the maximal ideal (i.e. if $C$ is Buchsbaum),
then the necessary and sufficient condition for the existence of such a curve
is given in terms of the dimensions of the components of $M$
(Proposition~\ref{degree-Buchsbaum}).
Similarly, if $\cal L$ contains any curve
with maximal corank then the minimal curve has
maximal corank as well, since any
basic double link increases $r_a $ by exactly 1 and $e$ by at least 1
(\cite[Lemma~1.14]{BM1}).  Hence again it reduces to a question of
the dimensions of the module components.

It is also true that $\cal L$ contains curves of maximal rank
(i.e.\ $r_0 < \alpha$) if and only if the minimal curve in
$\cal L$ has maximal rank, since a
basic double link increases $r_o$ by exactly 1 and increases $\alpha$ by at
most 1.  (This was first observed in \cite[Theorem~2.1]{BM3}.)  Then it follows
from Theorem~\ref{coh-prop} that if the deficiency module $M$ associated to
$\cal L$ has diameter 3 or more, and if $\cal L$ contains
any curve of maximal rank,
then $\cal L$ does not contain any curve achieving our bound on $\omega$.

Finally, we observe that it follows from Proposition~\ref{equal-coh} and
Corollary~\ref{almostBuchs-deg}
that if $\diam K < \diam M$ then $\cal L$ contains no curve
achieving our bound.  This immediately rules out a huge number of even liaison
classes, since ``most'' classes will have a module containing at least one pair
of consecutive components, for which multiplication by a general linear form is
injective.
\end{remark}

By using the results in Section~\ref{integral}, we can prove
similar results for the existence of integral curves with generators
of maximal degree.  However, because the theorems in that
section only go one direction, and because deformations do not in general
preserve integrality or degrees of generators, we can only
get necessity.

\begin{thm}  Suppose ${\cal L}$ is a liaison class whose minimal
curve $C_0$ satisfies $t(C_0) < t_1(C_0)$.  If ${\cal L}$ contains
an integral curve with a generator of maximal degree, then the
minimal curve $C_0$ has Hilbert function which satisfies
$$
\Delta^2 H(C_0, r_a(C_0) + 2) < \cdots < \Delta^2 H(C_0, \sigma(C_0) - 1) < 0.
$$
\end{thm}

\begin{proof} This follows immediately from Proposition~\ref{theta-coh}.
\end{proof}

As before, with extra assumptions on the liaison class, we
can give a cohomological criterion.

\begin{prop} Suppose ${\cal L}$ is a liaison class whose minimal
curve $C_0$ satisfies $t(C_0) < t_1(C_0)$.
\begin{enumerate}
\item[{\rm (}a.{\rm )}]  If $C_0$ has maximal corank, and if ${\cal L}$
contains integral curves with generators of maximal degree,
then
$$
h^1(\pthree, {\cal I}_{C_0}(t)) > 3 h^1(\pthree, {\cal I}_{C_0}(t+1))
		- 3 h^1(\pthree, {\cal I}_{C_0}(t+2))
		+   h^1(\pthree, {\cal I}_{C_0}(t+3)),
$$
for $t = r_a(C_0), \dots, r_o(C_0)$.
\item[{\rm (}b.{\rm )}]  If ${\cal L}_{n_1\dots n_r}$ is a Buchsbaum
liaison class containing integral curves with generators of maximal
degree, then $n_i > 3n_{i+1}$ for $i = 1, \dots, n_r$.
\end{enumerate}
\end{prop}

\begin{proof}  This follows exactly as before, using
Proposition~\ref{theta-coh}.
\end{proof}

\end{document}